\documentclass[iccmp]{ipbook}

\RequirePackage{amsthm,amsmath,amssymb,amsfonts}
\RequirePackage{hyperref}
\hypersetup{
    colorlinks,
    linktocpage,
    linkcolor={blue}
}
\RequirePackage{graphicx}
\RequirePackage{mathrsfs}
\RequirePackage{bm}
\RequirePackage{tabularx}
\RequirePackage{tikz}
\usetikzlibrary{arrows.meta,positioning,shapes.geometric,calc}

\startlocaldefs
\theoremstyle{plain}
\newtheorem{thm}{Theorem}[section]
\newtheorem{lem}[thm]{Lemma}

\theoremstyle{definition}
\newtheorem{defn}[thm]{Definition}

\newcommand{\veps}{\varepsilon}
\newcommand{\ri}{\mathrm{i}}
\newcommand{\og}{\omega}
\newcommand{\FGA}{\mathrm{FGA}}
\newcommand{\R}{\mathbb{R}}
\newcommand{\TT}{\mathbb{T}}
\newcommand{\C}{\mathbb{C}}

\newcommand{\bd}[1]{\boldsymbol{#1}}

\numberwithin{equation}{section}

\endlocaldefs

\firstpage{1}
\lastpage{30}

\title[Bridging scales by asymptotic analysis]%
{Bridging Scales: Asymptotic Analysis and AI-Assisted Formalization}

\author{Xu Yang}
\address{Department of Mathematics, University of California,
Santa Barbara, CA 93106, USA}
\email{xuyang@math.ucsb.edu}
\thanks{This paper is based on an invited plenary lecture delivered by the author at the International Congress of Chinese Mathematicians (ICCM 2025), January 2--8, 2026.}

\begin{document}

\begin{abstract}
Asymptotic analysis is one of the classical tools for bridging models across scales. Behind many such derivations lies a common
symbolic structure: an ansatz, a substitution, an order-by-order matching procedure, and the extraction of effective equations or interface conditions. This article revisits that structure through two representative bridges: the kinetic-to-fluid limit, illustrated by radiative transfer with interface layers and by neural-network approximations of Boltzmann equations, and the quantum-to-classical limit, illustrated by the Frozen Gaussian Approximation and its Dirac extension. We then explain why such derivations are natural candidates for AI-assisted formalization: their recurring symbolic structures can be organized, verified, and reused. In this sense, a carefully structured expository paper may serve not only as a review, but also as a mathematical seed for future AI-assisted
environments.
\end{abstract}

\maketitle

\section{Introduction}\label{sec:intro}
This paper revisits several classical asymptotic derivations from a
common perspective: they are not only bridges between models at
different physical scales, but also structured symbolic arguments.
This perspective guides both the review of classical multiscale
limits and the later discussion of AI-assisted formalization.

\subsection{Asymptotic analysis as a unifying lens across scales}
\label{subsec:intro-lens}
Much of applied mathematics is concerned with the passage between descriptions of the same physical system at different scales. A gas may be described by the positions and velocities of its molecules, by a kinetic density governed by a Boltzmann equation, or by the macroscopic density, velocity and temperature fields of fluid dynamics. A quantum particle may be described by a wave function solving the Schr\"odinger or Dirac equation, or, in the regime where the wavelength is small compared with the scale of variation of the potential, by classical trajectories carrying amplitudes and phases. In each case there is a small parameter (e.g., the Knudsen number, the mean free path, the rescaled Planck constant) whose vanishing marks the transition from a finer to a coarser model.

Asymptotic analysis is the mathematical tool that makes such transitions precise. Rather than solving the fine-scale model and then observing that its solutions resemble those of the coarse-scale model, one introduces an expansion of the solution in powers of the small parameter, substitutes it into the governing equation, and collects terms order by order. The leading order typically recovers the coarse-scale model; higher orders produce systematic corrections, and, when the geometry is nontrivial, boundary and interface layers that the coarse-scale model cannot see on its own. Classical examples include the Hilbert and Chapman--Enskog expansions connecting kinetic theory to hydrodynamics, the WKB method connecting wave mechanics to classical mechanics, and the matched asymptotic expansions of boundary-layer theory. It remains central to multiscale modeling today.

What unifies these examples is not the physics but the \emph{method}. The same handful of structural moves recurs: substitution of an ansatz, matching by orders of a small parameter, the imposition of solvability (Fredholm) conditions, the derivation of transport or diffusion equations for the surviving unknowns, and the resolution of layers by stretched variables. This observation is the starting point for the perspective of this paper.

\subsection{The structural pattern of asymptotic derivation}
\label{subsec:intro-pattern}
Viewed at the level of mathematical structure, many asymptotic derivations follow the same recipe. One first chooses an \emph{ansatz}, a parameterized form for the solution, such as a power series $f = f^{(0)} + \veps f^{(1)} + \veps^2 f^{(2)} + \cdots$, a single-phase WKB form $\og\, e^{\ri S/\veps}$, or a phase-space integral over Gaussian wave packets. One then \emph{substitutes} the ansatz into the governing equation and expands the resulting expression in powers of~$\veps$, often after a Taylor expansion of the potential or the collision kernel about a moving center. The next step is to \emph{match} terms order by order, obtaining at each order either an algebraic constraint, such as an eigenvalue relation or a solvability condition, or an evolution equation for one of the unknown profiles. Finally, one \emph{assembles} the resulting hierarchy into a closed coarse-scale model, sometimes supplemented by layer equations near boundaries or interfaces.

The individual steps within this recipe are drawn from a small, recurring vocabulary: power matching in~$\veps$; integration by parts under an oscillatory integral; the method of stationary phase; moment closure against a Maxwellian or against an eigenbasis; the Fredholm alternative; and changes of variables along characteristic flows. A reader who has carried out one such derivation in detail will recognize every one of these moves in the next, even when the physical setting is entirely different. The vocabulary is shared; only the nouns change.

\begin{figure}[ht]
  \centering
  \begin{tikzpicture}[
      node distance=5mm and 8mm,
      every node/.style={font=\footnotesize},
      box/.style={draw,rounded corners=1pt,align=center,
                  minimum height=8mm,inner sep=3pt,
                  text width=33mm},
      libbox/.style={draw,rounded corners=1pt,align=center,
                  minimum height=8mm,inner sep=3pt,
                  text width=33mm,fill=gray!10},
      arr/.style={-{Latex[length=2mm]},thick},
      sidearr/.style={-{Latex[length=2mm]},dashed,thick,gray!70}
  ]
    \node[box] (ansatz) {Ansatz \\ {\scriptsize (parameterized form)}};
    \node[box,below=of ansatz] (sub) {Substitute \\ {\scriptsize into governing equation}};
    \node[box,below=of sub] (match) {Match orders in~$\veps$ \\ {\scriptsize (solvability)}};
    \node[box,below=of match] (extract) {Effective models \\ {\scriptsize (interface conditions)}};
 
    \draw[arr] (ansatz) -- (sub);
    \draw[arr] (sub) -- (match);
    \draw[arr] (match) -- (extract);
 
    \node[libbox,right=22mm of ansatz] (libA) {Certified ansatz \\ templates};
    \node[libbox,right=22mm of sub] (libS) {Differentiation  \\ Taylor expansion};
    \node[libbox,right=22mm of match] (libM) {Order matching \\ solvability conditions};
    \node[libbox,right=22mm of extract] (libE) {Differential equations \\ matching diagrams};
 
    \draw[sidearr] (ansatz) -- (libA);
    \draw[sidearr] (sub) -- (libS);
    \draw[sidearr] (match) -- (libM);
    \draw[sidearr] (extract) -- (libE);
 
    \node[draw=none,right=4mm of libS,rotate=-90,anchor=south,
          inner sep=2pt,font=\scriptsize] (lib)
      {Reusable formal library};
    \draw[gray!70,thick]
      ($(libA.north east)+(2mm,0)$) --
      ($(libA.north east)+(5mm,0)$) --
      ($(libE.south east)+(5mm,0)$) --
      ($(libE.south east)+(2mm,0)$);
 
    \node[above=1mm of ansatz,font=\scriptsize\bfseries] {Classical asymptotic derivation};
    \node[above=1mm of libA,font=\scriptsize\bfseries] {AI-assisted formalization};
  \end{tikzpicture}
  \caption{An illustrative correspondence between the structural pattern of asymptotic analysis (left) and a possible AI-assisted view (right). The right-hand column is one tentative mapping among several; it is meant to make the connection between the classical
  recipe and machine-assisted formalization visible, not to define the components of any specific library. See Section~\ref{sec:ai} for the underlying discussion.}
  \label{fig:flowchart}
\end{figure}
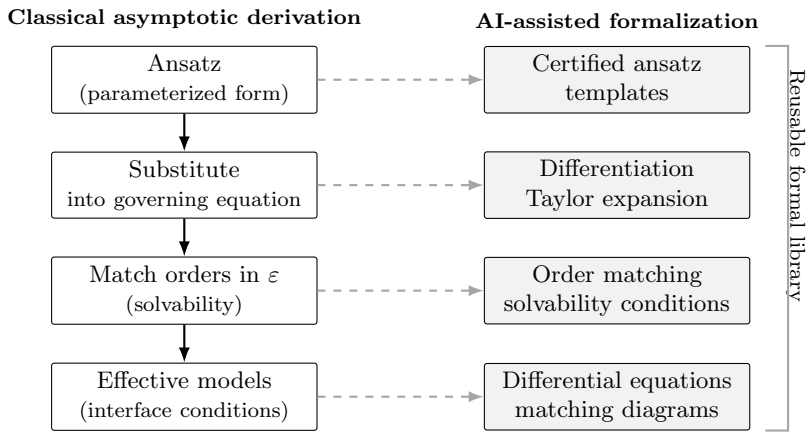

Figure~\ref{fig:flowchart} sketches this structural pattern as a flowchart. The classical recipe (ansatz, substitution, matching, extraction of effective equations) is shown on the left. The right-hand side is more speculative: it illustrates, by way of one possible correspondence, how the same symbolic moves might become targets for machine-assisted formalization, with each step suggesting candidate components for a growing certified library. The precise shape of such a library is left open here; the diagram is meant only to make the connection visible.
 
\subsection{The AI connection}
\label{subsec:intro-ai}
This vocabulary has a striking property: it is systematic in principle but error-prone in practice. Tracking which terms are $O(\veps)$ and which are $O(\veps^2)$, keeping signs and indices consistent through pages of oscillatory-integral manipulation, and verifying that the side conditions needed for integration by parts actually hold are often not where the main conceptual difficulty lies, but they are among the places where hand calculations are most vulnerable to error. This is precisely the kind of symbolic and structural work for which modern AI and formal verification tools are increasingly relevant: large language models can help draft the high-level structure of such arguments, while formal proof assistants such as Lean~\cite{Lean4} can, in principle, certify individual steps.

The perspective taken here is that an asymptotic derivation, viewed as a structured symbolic object built from a small recurring vocabulary, is a natural candidate for machine-assisted mathematics. Our aim is to make the case for this perspective, rather than to carry out such a formalization. Through two representative derivations and a closing discussion of one concrete framework for machine-assisted proof refinement, we suggest that the symbolic structure of asymptotic analysis is regular enough, and its recurring moves sufficiently few, that certified and reusable libraries for asymptotic derivations may be a productive direction. As a tentative organizing principle, we propose a \emph{reusable-formalization principle}: once a difficult member of a family of derivations has been formalized, the marginal cost of formalizing related derivations should drop substantially, because the same formal steps recur across the derivations.

It may help to separate what this article establishes from what it suggests. What we treat as established is structural: the derivations reviewed here are assembled from the same small vocabulary of symbolic moves, made explicit in each case. What we suggest is the possible formalization payoff: this vocabulary may be encoded as certified, reusable primitives in a proof assistant, and certifying one difficult derivation may lower the cost of formalizing related ones. The first is a claim about the mathematics, supported by the examples below; the second is a proposed direction for machine-assisted formalization, motivated by the repeated structure of these examples.

We develop this perspective through two classical bridges. The first, from kinetic to fluid, is reviewed through the diffusion limit of a radiative transfer equation with interfaces, and then connected to a recent error analysis of physics-informed neural networks for Boltzmann equations, a setting in which the asymptotic structure of the classical derivation reappears in the analysis of an AI-based solver. The second, from quantum to continuum, is developed in depth through the Frozen Gaussian Approximation, including its extension to the Dirac equation and the asymptotic preservation of the non-relativistic limit. The derivation of the Frozen Gaussian amplitude equation is the technical centerpiece of the paper and the concrete object to which the formalization discussion returns.

There is also a broader motivation for writing this paper in the present form. As AI-assisted mathematical systems become more powerful, research articles may increasingly serve a dual role. They will remain expository and archival documents for human readers, while also potentially serving as structured inputs to AI environments that verify, recombine, and extend existing arguments. A carefully written derivation, with its ansatz, matching procedure, solvability conditions, side assumptions, and reusable lemmas made explicit, can function as a mathematical seed: it records not only conclusions, but also reusable structures of reasoning that future machine-assisted systems may formalize, verify, and extend. In this sense, a review article can do more than summarize known results; it can identify symbolic patterns that may later be preserved, certified, and reused. This perspective motivates the way the material is presented below.

\smallskip

\noindent {\bf Outline.} Section~\ref{sec:kinetic-fluid} reviews the kinetic-to-fluid passage, including the diffusion limit of the radiative transfer equation, interface conditions obtained through an interface-layer analysis, and recent error and asymptotic-preservation estimates for physics-informed neural networks approximating Boltzmann equations. The goal is to show how the same asymptotic structure appears in both classical derivations and learning-based numerical methods. Section~\ref{sec:fga} develops the quantum-to-continuum passage through the Frozen Gaussian Approximation, presenting its formulation, the derivation of the amplitude equation, the extension to the Dirac equation, and the non-relativistic limit. This section serves as the main technical example. Section~\ref{sec:ai} turns to machine-assisted formalization: it describes the angelic-execution framework of~\cite{MathEye} and returns to the Frozen Gaussian derivation as a concrete candidate for the kind of structured argument that could be made explicit, certified, and reused. Section~\ref{sec:conclusion} concludes.

Throughout, the emphasis is on the \emph{structure} of the derivations rather than on the sharpest available hypotheses. Results are stated in the form most useful for the present discussion, and the original papers are cited for more complete assumptions and proofs.

\section{From kinetic to fluid}\label{sec:kinetic-fluid}
The passage from a kinetic description, in which a density is resolved in both physical space and velocity (or direction), to a fluid description in which only macroscopic moments survive, is the oldest and best-understood instance of the asymptotic bridges considered in this paper. We review it in a concrete setting (the radiative transfer equation with an interface) where the derivation produces not only a diffusion equation but also non-trivial transmission conditions at the interface, computed through an interface-layer analysis. We then turn to a modern counterpart, in which the same asymptotic structure governs the error analysis of a neural-network solver.

\subsection{The transport equation and its diffusive scaling}
\label{subsec:transport}
Consider the propagation of the phase-space energy density of waves in a heterogeneous medium with weak random fluctuations in the high frequency regime. In a planar-symmetric medium the three-dimensional radiative transfer equation reduces to a one-dimensional transport equation for the density $a(t,x,\mu)$, where $x\in\R$ is position and $\mu\in[-1,1]$ is the cosine of the angle between the direction of propagation and the~$x$-axis (e.g., \cite{yang2006numerical, jin2008domain, JLY2008}). Writing $v(x)$ for the background sound speed and $\sigma(x,\mu,\mu')$ for the differential scattering cross section, and rescaling time and space by the scattering mean free path~$\veps$, the equation takes the diffusive scaling
\begin{equation}\label{eq:rte}
  \veps^2\,\partial_t a
  + \veps\,\mu\,\partial_x a
  = \int_{-1}^{1}\sigma(x,\mu,\mu')
      \bigl[a(t,x,\mu') - a(t,x,\mu)\bigr]\,d\mu' .
\end{equation}
The total scattering cross section is $\Sigma(x) = \int_{-1}^1 \sigma(x,\mu,\mu')\,d\mu'$, and the right-hand side is the scattering operator, which relaxes $a$ toward its angular average. The factor $\veps^2$ in front of the time derivative and $\veps$ in front of the transport term encode the diffusive scaling: many scattering events occur on the time scale of interest, and the question is what effective equation governs the angular average of $a$ as $\veps\to0$.

\subsection{Hilbert expansion and the diffusion limit}
\label{subsec:hilbert}
The ansatz is the Hilbert expansion
\begin{equation}\label{eq:hilbert}
  a(t,x,\mu) = a^{(0)}(t,x,\mu) + \veps\,a^{(1)}(t,x,\mu)
             + \veps^2\,a^{(2)}(t,x,\mu) + \cdots .
\end{equation}
Substituting~\eqref{eq:hilbert} into~\eqref{eq:rte} and matching powers of~$\veps$ realizes the structural pattern of Section~\ref{subsec:intro-pattern} in its purest form. At order $\veps^0$ the scattering operator annihilates~$a^{(0)}$, which forces $a^{(0)}=a^{(0)}(t,x)$ to be independent of the angle~$\mu$: the leading-order density is isotropic. At order $\veps^1$ one solves an auxiliary equation for the first correction,
\begin{equation}\label{eq:a1}
  a^{(1)}(t,x,\mu) = -\frac{\theta(\mu)}{\Sigma}\,
     \partial_x a^{(0)}(t,x),
  \qquad
  \theta(\mu) - \int_{-1}^{1}\frac{\sigma(\mu,\mu')}{\Sigma}\,
     \theta(\mu')\,d\mu' = \mu ,
\end{equation}
the auxiliary function $\theta$ encoding the angular response of the medium. At order $\veps^2$ a solvability (Fredholm) condition, the requirement that the equation for $a^{(2)}$ have a solution, closes the hierarchy and produces a diffusion equation for the isotropic leading order. We record the outcome.

\begin{thm}[Diffusion limit \cite{JLY2008}]\label{thm:diffusion-limit}
Let $a$ solve~\eqref{eq:rte} with scattering cross section~$\sigma$ and total cross section~$\Sigma$ bounded and bounded below, and with well-prepared initial data. Then, away from boundaries and interfaces, the leading-order density $a^{(0)}(t,x)$ in the expansion~\eqref{eq:hilbert} is independent of~$\mu$ and satisfies the diffusion equation
\begin{equation}\label{eq:diffusion}
  2\,\partial_t a^{(0)}
  - \partial_x\!\left(D\,\partial_x a^{(0)}\right) = 0,
  \qquad
  D = \frac{1}{\Sigma}\int_{-1}^{1}\mu\,\theta(\mu)\,d\mu ,
\end{equation}
with $\theta$ the solution of the auxiliary problem~\eqref{eq:a1}.
\end{thm}

The diffusion constant~$D$ is not an input but an output of the derivation: it is computed from the microscopic scattering data through the auxiliary function~$\theta$. This is characteristic of asymptotic bridges: the coarse model inherits effective coefficients determined by the fine model.

\subsection{Interface conditions via interface-layer analysis}
\label{subsec:interface}
The diffusion limit of Theorem~\ref{thm:diffusion-limit} holds in the interior of a homogeneous medium. The situation of real interest is when the sound speed~$v(x)$ is discontinuous across an interface, say at $x=0$, separating two media. Waves crossing the interface undergo reflection and transmission, and a naive matching of the two interior diffusion equations across $x=0$ is incorrect: an interface layer forms, in which the density is genuinely anisotropic, and the correct conditions linking the two diffusion solutions must be extracted from a boundary-layer analysis of the kinetic equation~\cite{JLY2008,BalRyzhik}.

At the kinetic level the two solutions are coupled by reflection and transmission operators that redistribute energy among directions. The angle of transmission is fixed by Snell's law,
\begin{equation}\label{eq:snell}
  \frac{\sqrt{1-\mu_1^2}}{v_1} = \frac{\sqrt{1-\mu_2^2}}{v_2},
\end{equation}
relating the cosine~$\mu_1$ of the incident angle to the cosine~$\mu_2$ of the transmitted angle, while reflection is specular. The response of each half-space to a prescribed incoming density is encoded in operators built from the Chandrasekhar $H$-function~\cite{Chandrasekhar}, the classical special function of radiative transfer. Resolving the interface layer reduces to an eigenvalue problem posed at $x=0$, whose solution yields a single effective transmission constant.

\begin{thm}[Effective interface condition \cite{JLY2008,BalRyzhik}]
\label{thm:interface}
Across a discontinuity of the sound speed at $x=0$, the leading-order diffusion densities $u_1 = a^{(0)}|_{x>0}$ and $u_2 = a^{(0)}|_{x<0}$ satisfy the matching conditions 
\begin{equation}\label{eq:interface-jump}
  u_2(0^-) = \alpha\,u_1(0^+),
  \qquad
  D_2\,\partial_x u_2(0^-) = D_1\,\partial_x u_1(0^+),
\end{equation}
where the flux is continuous and the density jumps by a factor~$\alpha$. For specular reflection with regular transmission, $\alpha = v_1^2/v_2^2$; for diffuse transmission, $\alpha$ is determined by an eigenvalue problem at the interface expressed through the Chandrasekhar $H$-function and the reflection and transmission operators.
\end{thm}

The derivation of Theorem~\ref{thm:interface} is the part of the kinetic-to-fluid story that most resembles the quantum case treated later: it is not a matter of formally setting $\veps=0$, but of constructing a layer correction whose matching to the interior solution produces a non-trivial effective condition. The two interface conditions (regular and diffuse) give visibly different interior profiles near $x=0$; the diffuse case exhibits a genuine interface layer that the regular case does not, as confirmed by the numerical comparison in~\cite{JLY2008}.

\subsection{Physics-informed neural networks for Boltzmann equations}
\label{subsec:pinn}
The classical derivation above takes a kinetic model and produces a fluid model. A complementary modern question is whether a \emph{numerical} solver based on machine learning respects the same asymptotic structure. We summarize recent results of Abdo, Chai, Hu and the author~\cite{ACHY} on physics-informed neural networks (PINNs) for the Boltzmann equation, which show that the asymptotic analysis of the classical bridge reappears, almost unchanged, in the error analysis of the AI-based method.

The classical Boltzmann equation $\partial_t f + \xi\cdot\nabla_x f = Q(f,f)$ governs the density $f(t,x,\xi)$ of particles with velocity~$\xi$. Linearizing about a global Maxwellian $M$ through $f = M + M^{1/2}u$ yields the linearized equation
\begin{equation}\label{eq:linboltz}
  \partial_t u + \xi\cdot\nabla_x u + \nu(\xi)\,u
  = K(u) + \Gamma(u,u),
\end{equation}
with a collision frequency~$\nu$ satisfying $\nu_0(1+|\xi|)^\gamma\le\nu(\xi)\le\nu_1(1+|\xi|)^\gamma$ for some $\gamma\in[0,1]$, a bounded nonlocal operator~$K$, and a quadratic term~$\Gamma$. A PINN approximates the solution by a neural network~$u_\theta$ whose parameters~$\theta$ are trained to make a \emph{residual} small. On the unbounded velocity domain one truncates to a cube~$Q_R$ of diameter~$R$ and defines the PDE residual
\begin{equation}\label{eq:residual}
  \mathcal{R}_i[R;\theta]
  = \partial_t u_\theta + \xi\cdot\nabla_x u_\theta + \nu(\xi)\,u_\theta
    - K(u_\theta\chi_{Q_R})
    - \Gamma(u_\theta\chi_{Q_R}, u_\theta\chi_{Q_R}),
\end{equation}
together with initial and (periodic) boundary residuals. The square of the \emph{generalization error} $\mathcal{E}_G[R;\theta]$ is the sum of the squared continuous $L^2$ residual norms on the corresponding space--time, initial, and boundary domains, together with the penalty $\lambda_R(1+R)^\gamma\int_0^T\|u_\theta\|_{H^2(\TT^3\times Q_R)}^2\,dt$, where $\lambda_R>0$. The \emph{total error} $\mathcal{E}[R;\theta]$ is the $L^2$ distance between the network and the true solution on $[0,T]\times\TT^3\times Q_R$. The central question is whether the continuous residuals control this solution error, accounting for velocity truncation and for the sampled loss used in training.

\begin{thm}[Total error controlled by residual \cite{ACHY}]
\label{thm:pinn-total}
Under the collision-kernel assumptions of~\cite{ACHY}, let $u$ solve~\eqref{eq:linboltz} on $[0,T]\times\TT^3\times\R^3$ with $M_u := \|(1+|\xi|)^{\gamma/2}u\|_{L^\infty(0,T;L^\infty(\TT^3;L^2(\R^3)))}$ finite, let $u_\theta$ be a $\tanh$ neural network, and define the velocity tail
\[
  \tau(R)^2 = \int_0^T\!\!\int_{\TT^3}\!\int_{\R^3\setminus Q_R}(1+|\xi|)^\gamma|u|^2\,d\xi\,dx\,dt ,
\]
so that $\tau(R)\to0$ as $R\to\infty$. For $R\ge1$ and $\lambda_R>0$, there is a constant~$C$, depending on~$T$, $M_u$, and the fixed collision-kernel bounds, but independent of~$R$, $\theta$, and~$\lambda_R$, such that
\begin{equation}\label{eq:pinn-bound}
  \mathcal{E}[R;\theta]^2
  \le C\exp\!\left(C\,\frac{\mathcal{E}_G[R;\theta]^2}{\lambda_R}\right)
      \bigl(\mathcal{E}_G[R;\theta]^2 + \tau(R)^2\bigr).
\end{equation}
\end{thm}

The tail term $\tau(R)$ accounts for the velocities omitted by the truncation, which the residual on~$Q_R$ does not see. Estimate~\eqref{eq:pinn-bound} follows by retaining this term in the energy estimate of~\cite[Theorem~4.1]{ACHY}, rather than absorbing it into the residual through a suitable choice of~$R$. We also retain the dependence of the Gronwall factor on the regularization parameter. Consequently, convergence along a sequence of networks and truncation radii is guaranteed when the full right-hand side of~\eqref{eq:pinn-bound} tends to zero; $\mathcal{E}_G\to0$ and $R\to\infty$ alone do not control that factor.

To relate this estimate to training, let $\mathcal{L}_N[R;\theta]$ denote a consistently weighted quadrature approximation of $\mathcal{E}_G[R;\theta]^2$, including the penalty. A quadrature estimate takes the form
\[
  \bigl|\mathcal{E}_G[R;\theta]^2-\mathcal{L}_N[R;\theta]\bigr|
  \le \delta_N(R,\theta).
\]
For fixed~$R$ and a fixed network, the midpoint rule gives $\delta_N(R,\theta)\to0$ under the residual regularity assumptions of~\cite[Remark~2]{ACHY}. Its error constants depend on derivatives of the residuals. When~$R$ and~$\theta$ also vary, the sampling must therefore be refined so that this error vanishes along the same sequence. A small sampled loss controls the continuous generalization error only up to this quadrature error.

Under the stronger weighted regularity assumption $(1+|\xi|)^\gamma u\in H^4([0,T]\times\TT^3\times\R^3)$, the approximation result in~\cite[Theorem~3.2]{ACHY} yields arbitrarily small generalization error with suitable choices of truncation radius, regularization parameter, and network size. This existence result does not itself control the Gronwall factor in~\eqref{eq:pinn-bound}. The proof of the residual estimate rests on local estimates for the nonlocal operators~$K$ and~$\Gamma$ on the truncated domain, the analytic analogue of the moment estimates that close the Hilbert hierarchy in Section~\ref{subsec:hilbert}.

\subsection{Micro--macro decomposition and asymptotic preservation}
\label{subsec:ap}

The deepest point of contact between the classical and AI-based analyses is the multiscale regime. Consider the rescaled linear Boltzmann equation
\begin{equation}\label{eq:multiscale}
  \veps\,\partial_t f^\veps + \xi\cdot\nabla_x f^\veps
  = \frac{1}{\veps}\,L(f^\veps),
\end{equation}
where the linear collision operator is
\begin{equation}\label{eq:multiscale-L}
  L(f) = \int_{\R^3}\alpha(\xi,\xi')
     \bigl(M(\xi)f(\xi') - M(\xi')f(\xi)\bigr)\,d\xi',
\end{equation}
with a symmetric, bounded, strictly positive kernel~$\alpha$. As $\veps\to0$ this equation has a diffusion limit, exactly as the radiative transfer equation did. A direct PINN discretization of \eqref{eq:multiscale} degrades as $\veps\to0$, because the stiff factor~$1/\veps$ amplifies residual errors. The remedy, following the asymptotic-preserving philosophy~\cite{jin2023asymptotic}, is to build the micro--macro decomposition into the network. One writes
\begin{equation}\label{eq:micromacro}
  f^\veps(t,x,\xi) = \rho^\veps(t,x)\,M(\xi) + \veps\,g^\veps(t,x,\xi),
  \qquad
  \int_{\R^3} g^\veps\,d\xi \equiv 0 ,
\end{equation}
separating the equilibrium part $\rho^\veps M$ from the microscopic fluctuation~$g^\veps$. Substituting~\eqref{eq:micromacro} into~\eqref{eq:multiscale} and projecting yields a coupled macro--micro system, and as $\veps\to0$ the macroscopic density~$\rho^\veps$ converges to the solution~$\rho^0$ of a diffusion equation. Parameterizing $\rho^\veps$ and~$g^\veps$ separately by neural networks (rather than parameterizing $f^\veps$ directly) yields an asymptotic-preserving PINN whose accuracy survives the limit.

\begin{thm}[Asymptotic preservation \cite{ACHY}]
\label{thm:ap}
For the micro--macro PINN associated with~\eqref{eq:multiscale}, the total error of the macroscopic density~$\rho^\veps$ is bounded by a multiple of the generalization error~$\mathcal{E}_G(\mathcal{F}^\veps)$ of the coupled system, uniformly in~$\veps$. Consequently, since $\|\rho^\veps - \rho^0\| = O(\veps)$, minimizing the generalization error yields a network that accurately approximates the diffusion limit~$\rho^0$ as $\veps\to0$.
\end{thm}

Theorem~\ref{thm:ap} is the statement that the diagram of approximations commutes: taking the network limit and taking the fluid limit can be interchanged. This is the same commuting-diagram structure that, in Section~\ref{sec:fga}, will express the compatibility of the Frozen Gaussian Approximation with the non-relativistic limit. The fact that it appears here, in the analysis of a neural-network solver, illustrates the perspective taken here: the asymptotic structure of the classical derivation is not an artifact of the classical method but appears to be an invariant of the problem, governing the AI-based method as it governs the pencil-and-paper one.

\begin{figure}[ht]
  \centering
  \includegraphics[width=0.95\textwidth]{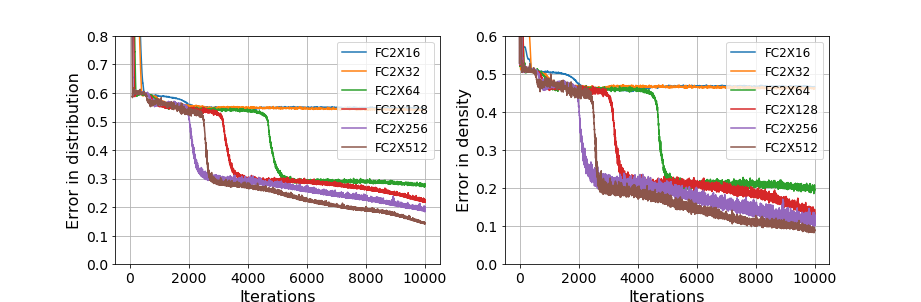}
  \caption{Convergence of the physics-informed neural network approximation for the Boltzmann equation as training proceeds. The left panel shows the error in the distribution function and the right panel shows the error in the density, plotted against the number of training iterations for fully-connected networks of fixed depth and varying width (16, 32, 64, 128, 256, 512 neurons per layer). As the width grows the asymptotic error decreases and the training loss stabilizes at progressively smaller values, qualitatively consistent with Theorem~\ref{thm:pinn-total}. Reproduced from~\cite{ACHY}.}
  \label{fig:pinn}
\end{figure}

\subsection{What the kinetic-to-fluid story illustrates}
\label{subsec:k2f-discussion}
Two derivations have appeared in this section. The classical one takes the radiative transfer equation, applies a Hilbert expansion, and produces a diffusion equation with $H$-function interface conditions. The modern one takes the linearized Boltzmann equation, applies a micro--macro decomposition, and produces error and asymptotic-preservation estimates for a neural-network solver. Superficially these are different enterprises, one analytic and one computational. Structurally they are the same: both rest on separating an equilibrium (isotropic, Maxwellian) component from a fluctuation, on a solvability condition that closes the equation for the equilibrium part, and on the commutativity of the coarse-graining limit with the construction at hand. The vocabulary of Section~\ref{subsec:intro-pattern} (moment closure, Fredholm solvability, commuting limits) is common to both. We now turn to a setting where the derivation is substantially more intricate, and where the case for AI assistance is correspondingly stronger.

\section{From quantum to classical}\label{sec:fga}
We now develop the second asymptotic bridge of the paper, from quantum wave mechanics to classical phase-space dynamics, in substantially more depth. The central construction is the Frozen Gaussian Approximation (FGA); we derive it for the Dirac equation and establish the asymptotic preservation of the non-relativistic limit. The derivation of the amplitude equation, carried out in Section~\ref{subsec:fga-derivation}, is the technical centerpiece of the paper and the object to which the discussion of formalization in Section~\ref{sec:ai} returns.

\subsection{The semiclassical Schr\"odinger equation}
\label{subsec:schrodinger}
The simplest quantum-to-classical bridge concerns the semiclassical Schr\"odinger equation
\begin{equation}\label{eq:schrodinger}
  \ri\veps\,\partial_t\psi^\veps
  = -\frac{\veps^2}{2}\Delta\psi^\veps + V(x)\,\psi^\veps,
  \qquad 0 < \veps \ll 1,
\end{equation}
where $\veps$ is the rescaled Planck constant and $V$ a smooth potential. As $\veps\to0$ the solution oscillates on the scale~$\veps$ in both space and time, which makes direct numerical simulation prohibitively expensive: resolving the oscillations requires a mesh finer than~$\veps$ in every dimension. The asymptotic methods of this section sidestep the cost by representing the oscillation analytically and evolving only the slowly varying data that ride on top of it.

\subsection{WKB and the breakdown at caustics}
\label{subsec:wkb}
The classical ansatz is the WKB (geometric optics) form, a single oscillatory wave with slowly varying amplitude and phase,
\begin{equation}\label{eq:wkb}
  \psi^\veps(t,x) \approx A(t,x)\,
    \exp\!\left(\frac{\ri}{\veps}\,S(t,x)\right).
\end{equation}
Substituting~\eqref{eq:wkb} into~\eqref{eq:schrodinger} and matching powers of~$\veps$ (the same move as in Section~\ref{subsec:hilbert}, now applied to an oscillatory rather than a power-series ansatz) gives at order~$\veps^0$ the eikonal equation
\begin{equation}\label{eq:eikonal}
  \partial_t S + \tfrac12|\nabla S|^2 + V = 0,
\end{equation}
a Hamilton--Jacobi equation, and at order~$\veps^1$ the transport equation
\begin{equation}\label{eq:transport-wkb}
  \partial_t A + \nabla S\cdot\nabla A
  + \tfrac12 A\,\Delta S = 0
\end{equation}
for the amplitude. The eikonal equation is solved by the method of characteristics: phase information is carried along classical trajectories $(Q(t),P(t))$ obeying Hamilton's equations $\dot Q = P$, $\dot P = -\nabla V(Q)$. The construction is exact and cheap, until two trajectories cross. At such a \emph{caustic} the map $q\mapsto Q(t,q)$ from initial to current position becomes singular, its Jacobian vanishes, the amplitude~$A$ in~\eqref{eq:transport-wkb} blows up, and the single-phase ansatz~\eqref{eq:wkb} fails. Caustics are generic (any focusing potential produces them) so the breakdown is not a pathology to be assumed away but the central difficulty the more sophisticated methods are designed to overcome.

\subsection{Gaussian beam method}
\label{subsec:gb}
The Gaussian beam method removes the caustic singularity by allowing the phase to be complex. Around each central trajectory one introduces a locally quadratic phase whose Hessian~$M(t)$ is a complex symmetric matrix with positive-definite imaginary part, so that the beam decays like a Gaussian away from its center. The Hessian evolves by a matrix Riccati equation derived, once again, by matching orders of~$\veps$ in a Taylor expansion of the potential about the central trajectory. Because $\operatorname{Im}M(t)$ remains positive definite for all time, the amplitude never blows up, and a superposition of Gaussian beams represents the solution accurately even through caustics~\cite{JWY2008,JWY2011,JWYH2010,WY2013,JMS}. The price is that an individual beam spreads as $\operatorname{Im}M(t)$ decreases, so that accuracy over long times requires either reinitialization or a large number of beams. This tension (between the analytic convenience of a moving Gaussian and the spreading of its width) is exactly what the Frozen Gaussian Approximation resolves.

\subsection{The Frozen Gaussian Approximation}
\label{subsec:fga-formulation}
The idea of the FGA, originating in the chemistry literature with Herman and Kluk~\cite{HK} and given a rigorous mathematical justification in \cite{swart2009mathematical, LuYang}, is to keep the Gaussian width \emph{frozen} (fixed once and for all) and to recover accuracy by superposing such frozen Gaussians over all of phase space, with a slowly varying amplitude that absorbs the spreading. The solution is represented as a phase-space integral over Gaussian wave packets centered on classical trajectories.

We state the construction for the scalar Schr\"odinger case~\eqref{eq:schrodinger}, posed in spatial dimension~$d$; the Dirac version used in the following subsections has the same structure~\cite{CLY}. In the scalar case~\cite{HK}
\begin{equation}\label{eq:fga-schr}
  \psi^\veps_{\FGA}(t,x)
  = \frac{1}{(2\pi\veps)^{3d/2}}\int_{\R^{3d}}
      a(t,q,p)\,
      \exp\!\left(\frac{\ri}{\veps}\,\Phi(t,x,y,q,p)\right)
      \psi^\veps_I(y)\,dy\,dq\,dp,
\end{equation}
with the FGA phase
\begin{equation}\label{eq:fga-phase}
  \begin{aligned}
  \Phi ={}& S(t,q,p)
       + \frac{\ri}{2}|x - Q(t,q,p)|^2
       + P(t,q,p)\cdot\bigl(x - Q(t,q,p)\bigr)\\
       &+ \frac{\ri}{2}|y - q|^2
       - p\cdot(y - q).
  \end{aligned}
\end{equation}
Here $(q,p)$ ranges over initial phase-space points, $(Q,P)$ are the classical trajectories issuing from them, and the integral over~$y$ projects the initial data onto the frozen Gaussian centered at~$q$. The Gaussian factors $\frac{\ri}{2}|x-Q|^2$ and $\frac{\ri}{2}|y-q|^2$ have \emph{fixed} width, in contrast to the evolving Hessian of a Gaussian beam; this is the sense in which the Gaussian is frozen.

\begin{defn}[FGA ansatz]\label{def:fga}
The FGA representation of the solution to a semiclassical evolution equation is a phase-space integral of the form~\eqref{eq:fga-schr}, in which the centers $(Q,P)$ evolve along the Hamiltonian flow of the classical symbol, the action~$S$ accumulates along the flow, and the amplitude~$a$ solves a transport equation derived by matching the ansatz against the equation to order~$\veps$.
\end{defn}

The centers and the action are governed by the Hamiltonian flow and its accompanying action integral,
\begin{equation}\label{eq:fga-flow}
  \dot Q = \partial_P h(Q,P), \quad
  \dot P = -\partial_Q h(Q,P), \quad
  \dot S = P\cdot\partial_t Q - h(Q,P),
\end{equation}
with $Q(0)=q$, $P(0)=p$, $S(0)=0$, where $h$ is the classical Hamiltonian (for Schr\"odinger, $h = |P|^2/2 + V(Q)$). What remains, and what carries all the analytic content, is the equation for the amplitude~$a$.

\subsection{Derivation of the FGA amplitude equation}
\label{subsec:fga-derivation}
This subsection is the technical heart of the paper. The derivation of the amplitude equation is where the structural moves of asymptotic analysis appear in their most intricate and most mechanizable form, and it is the concrete object to which Section~\ref{sec:ai} returns. We present it for the scalar case and indicate the Dirac modifications in Section~\ref{subsec:fga-dirac}; the full matrix-valued derivation is in~\cite{CLY}.

For the scalar analysis below, we assume that $V\in C^\infty(\R^d;\R)$ is \emph{subquadratic}, meaning that
\[
  \sup_{x\in\R^d}|\partial_x^\beta V(x)|<\infty
  \qquad\text{for every }|\beta|\ge2.
\]
On a fixed time interval $[0,T]$, this condition ensures that the Hamiltonian flow has uniformly bounded positive-order derivatives in $(q,p)$. The associated symplectic structure implies that $Z=\partial_z(Q+\ri P)$ is invertible and that $Z^{-1}$ and its phase-space derivatives are uniformly bounded~\cite{swart2009mathematical}. These bounds are the analytic hypotheses used in the oscillatory-integral estimates below.

One substitutes the ansatz, with a first-order correction~$b$,
\begin{equation}\label{eq:fga-ansatz}
  \psi^\veps(t,x)
  = \frac{1}{(2\pi\veps)^{3d/2}}\int_{\R^{3d}}
      \bigl(a(t,q,p) + \veps\, b(t,q,p)\bigr)\,
      e^{\ri\Phi/\veps}\,\psi^\veps_I(y)\,dy\,dq\,dp,
\end{equation}
into the evolution equation, and expands everything in powers of~$\veps$. The differentiations of the phase produce the key identities
\begin{equation}\label{eq:phase-derivs}
  \partial_t\Phi
  = \partial_t S - P\cdot\partial_t Q
    + (x - Q)\cdot\partial_t(P - \ri Q),
  \qquad
  \nabla_x\Phi = \ri(x - Q) + P,
\end{equation}
and the scalar potential is Taylor-expanded about the moving center~$Q$,
\begin{equation}\label{eq:V-taylor}
  V(x) = V(Q) + (x - Q)\cdot\partial_Q V(Q)
       + \tfrac12 (x - Q)^2 : \partial_Q^2 V(Q)
       + R_3(x,Q).
\end{equation}
Here the integral remainder is
\[
  \begin{aligned}
  R_3(x,Q)&=\sum_{|\alpha|=3}r_\alpha(x,Q)(x-Q)^\alpha,\\
  r_\alpha(x,Q)&=\frac{3}{\alpha!}\int_0^1(1-s)^2
       \partial^\alpha V\bigl(Q+s(x-Q)\bigr)\,ds.
  \end{aligned}
\]
The coefficients $r_\alpha(x,Q(t,q,p))$ and all their derivatives in $(x,q,p)$ are bounded uniformly on $[0,T]$. The matrix potential in the Dirac case admits the analogous entrywise expansion under the hypotheses of~\cite{CLY}.

The powers of $(x-Q)$ generated by~\eqref{eq:V-taylor} become small \emph{under the oscillatory integral}. Each integration by parts in phase space produces a factor of~$\veps$, but the derivative can also act on another factor of $(x-Q)$, reducing the polynomial degree by two in a single step. These contractions determine the order estimate in the following lemma. We use the complex differential operator $\partial_z=\partial_q-\ri\,\partial_p$ and the matrix $Z=\partial_z(Q+\ri P)$.

\begin{lem}[Oscillatory-integral identities; cf.~\cite{swart2009mathematical, LuYang,CLY}]
\label{lem:osc}
Let the phase~$\Phi$ be given by~\eqref{eq:fga-phase}, with the Hamiltonian flow and action~\eqref{eq:fga-flow}. Assume that the flow has bounded positive-order phase-space derivatives and that $Z^{-1}$ and its phase-space derivatives are bounded on $[0,T]\times\R^{2d}$, as ensured by the scalar subquadratic hypothesis above. For an amplitude~$f=f(t,x,q,p)$, define the operator
\[
  \mathcal{I}_\veps[f]\,\psi^\veps_I
  = \frac{1}{(2\pi\veps)^{3d/2}}\int_{\R^{3d}} f\,e^{\ri\Phi/\veps}\,\psi^\veps_I(y)\,dy\,dq\,dp .
\]
The integral is understood as an oscillatory integral, initially on Schwartz data and then by its bounded extension to~$L^2$. Write $f\sim g$ if $\mathcal{I}_\veps[f]=\mathcal{I}_\veps[g]$, and $f\sim O(\veps^k)$ if $\|\mathcal{I}_\veps[f]\|_{L^2\to L^2}\le C_T\veps^k$ for $0<\veps\le1$, uniformly for $t\in[0,T]$. Let $w$ and $G$ be smooth vector- and matrix-valued coefficients whose derivatives in $(x,q,p)$, including order zero, are uniformly bounded; if they depend on~$\veps$, these bounds are independent of~$\veps$. Coefficients independent of~$x$ are included in this class. All phase-space differentiations below keep~$x$ fixed. With Einstein's summation convention,
\begin{align}
  w\cdot(x - Q)
    &\sim -\veps\,\partial_{z_k}\!\bigl(w_j\,Z^{-1}_{jk}\bigr),
  \label{eq:osc-linear}\\
  (x - Q)\cdot G\,(x - Q)
    &\sim \veps\,\partial_{z_k}Q_l\,G_{lj}\,Z^{-1}_{jk}
      + O(\veps^2).
  \label{eq:osc-quadratic}
\end{align}
More generally, for a multi-index~$\alpha$ and any scalar coefficient~$c$ in this bounded-derivative class,
$c\,(x-Q)^\alpha \sim O(\veps^{\lceil|\alpha|/2\rceil})$. In particular, cubic and quartic terms are both $O(\veps^2)$ in this operator sense. The $O(\veps^2)$ remainder in~\eqref{eq:osc-quadratic} has the same meaning.
\end{lem}

Indeed, the Hamiltonian flow and action identities give $\ri\partial_z\Phi=Z(x-Q)$. Integration by parts yields~\eqref{eq:osc-linear}; applying it to the quadratic term produces the contraction in~\eqref{eq:osc-quadratic}. Iteration reduces the polynomial degree by one or two at each step, giving at least $\lceil|\alpha|/2\rceil$ factors of~$\veps$. The resulting coefficients remain in the bounded-derivative class and, viewed as symbols independent of~$y$, belong to the class $S[0;4d]$ of~\cite{swart2009mathematical}. The general $L^2$ bound in~\cite[Theorem~1, Eq.~(17)]{swart2009mathematical} allows symbols depending on $(x,y,q,p)$ and therefore gives the stated norm estimates. The permitted $x$-dependence includes the Taylor coefficients $r_\alpha$ above, so for a bounded-derivative amplitude~$a$ the remainder satisfies $aR_3\sim O(\veps^2)$.

The two identities of Lemma~\ref{lem:osc} are the workhorses: the first converts the linear term in the Taylor expansion~\eqref{eq:V-taylor} into an $O(\veps)$ contribution, and the second does the same for the quadratic term, each time at the cost of a holomorphic derivative and a factor of~$Z^{-1}$. With these in hand the matching proceeds in two stages. At order~$\veps^0$ one obtains the eikonal-type relation, which forces~$a$ to be an eigenvector of the principal symbol and fixes the action~$S$ to satisfy~\eqref{eq:fga-flow}. At order~$\veps^1$, after applying Lemma~\ref{lem:osc} to remove every factor of $(x-Q)$, the first-order correction~$b$ drops out (it is annihilated by the same eigenprojection that fixed~$a$), and what survives is a closed transport equation for~$a$.

\begin{thm}[FGA amplitude equation \cite{HK}]\label{thm:transport}
For the semiclassical Schr\"odinger equation~\eqref{eq:schrodinger}, the FGA amplitude in~\eqref{eq:fga-schr} satisfies the transport equation
\begin{equation}\label{eq:fga-amp-schr}
  \frac{d a}{dt}
  = \frac{a}{2}\,\operatorname{tr}\!\Bigl(Z^{-1}\bigl(\partial_z P
     - \ri\,\partial_z Q\,\partial_Q^2 V(Q)\bigr)\Bigr),
  \qquad a(0,q,p) = 2^{d/2},
\end{equation}
with $Z = \partial_z(Q + \ri P)$, $\partial_z = \partial_q - \ri\,\partial_p$, where, as in~\cite{LuYang}, the row index of a Jacobian is the derivative index: $(\partial_z Q)_{kj} = \partial_{z_k}Q_j$ and $Z_{kj} = \partial_{z_k}(Q_j + \ri P_j)$. The two contributions encode, respectively, the spreading of the phase-space volume element and the curvature of the potential along the flow. Since $\dot Z = \partial_z P - \ri\,\partial_z Q\,\partial_Q^2V(Q)$ by~\eqref{eq:fga-flow}, equivalently $a = (\det Z)^{1/2}$, with the branch of the square root chosen continuously in~$t$.
\end{thm}

Equation~\eqref{eq:fga-amp-schr} is the mechanism by which the Frozen Gaussian Approximation maintains accuracy despite the fixed Gaussian width: the amplitude~$a$ evolves so as to account for the spreading and interference of the underlying classical flow. The derivation just sketched consists of substituting the ansatz, expanding the phase, Taylor-expanding the potential, applying the oscillatory-integral identities to trade powers of $(x-Q)$ for powers of~$\veps$, projecting onto the relevant eigenspace, and identifying the remaining order-\(\veps\) terms as the transport equation for the amplitude. This is a representative instance of the structural pattern described in Section~\ref{subsec:intro-pattern}, and it involves enough sign, index, and order tracking that a hand derivation is genuinely error-prone. We return to this point in Section~\ref{sec:ai}.

\subsection{FGA for the Dirac equation}
\label{subsec:fga-dirac}
The Dirac equation describes a relativistic spin-$\tfrac12$ particle and, in the semiclassical regime, reads
\begin{equation}\label{eq:dirac}
  \ri\veps\,\partial_t\psi^\veps
  = \bigl(-\ri c\veps\,\widehat{\boldsymbol\alpha}\cdot\nabla
    - \widehat{\boldsymbol\alpha}\cdot A(x)
    + m\beta c^2 + V(x)\bigr)\psi^\veps,
\end{equation}
where $\psi^\veps\in\C^4$ is a four-spinor, $c$ is the speed of light, and $\widehat{\boldsymbol\alpha}$, $\beta$ are the Dirac matrices built from the $2\times2$ Pauli matrices~\cite{CLY}. The classical symbol is $D(q,p) = c\,\widehat{\boldsymbol\alpha}\cdot p + B(q)$, where
\[
  B(q)=-\widehat{\boldsymbol\alpha}\cdot A(q)+m\beta c^2+V(q)I_4
\]
and $I_4$ is the $4\times4$ identity matrix. The symbol~$D$ is Hermitian and has two \emph{double} eigenvalues
\begin{equation}\label{eq:dirac-eig}
  h_\pm(q,p) = \pm\lambda(q,p) + V(q), \qquad
  \lambda(q,p) = \sqrt{|pc - A(q)|^2 + m^2c^4},
\end{equation}
corresponding to the positive- and negative-energy branches. The new feature, absent in every previously analyzed case, is that each eigenvalue has multiplicity two: the corresponding eigenspace is two-dimensional, spanned by a pair of eigenvectors $r_{\pm,1}, r_{\pm,2}$. The amplitude is therefore not a scalar but a two-component vector $a = (a_1, a_2)^\top$ within each eigenspace, and the projection step in the derivation must resolve the amplitude within a degenerate eigenspace rather than along a single eigendirection.

The consequence is that the transport equation \eqref{eq:fga-amp-schr} is replaced by a coupled system
\begin{equation}\label{eq:dirac-amp}
  \frac{d}{dt}
  \begin{pmatrix} a_1 \\ a_2 \end{pmatrix}
  = \bigl(\mathcal{L} + \mathcal{M} + \mathcal{N}\bigr)
  \begin{pmatrix} a_1 \\ a_2 \end{pmatrix},
\end{equation}
in which $\mathcal{L}$, $\mathcal{M}$, $\mathcal{N}$ are explicit $2\times2$ matrices: $\mathcal{L}$ comes from the linear Taylor term via~\eqref{eq:osc-linear}, $\mathcal{M}$ from the quadratic term via~\eqref{eq:osc-quadratic}, and $\mathcal{N}$ from the variation of the eigenvectors $r_{\pm,j}$ along the flow. The structure of the derivation is identical to the scalar case; only the eigenprojection is more delicate, because it must keep track of how the two-dimensional eigenspace rotates as $(Q,P)$ evolve. This is the precise sense in which the Dirac case is the harder member of the family, and, by the reusable-formalization principle, the right one to formalize.

\subsection{Non-relativistic limit and asymptotic preservation}
\label{subsec:nonrel}
The Dirac equation carries, in addition to the semiclassical parameter~$\veps$, a second asymptotic parameter: the speed of light~$c$. As $c\to\infty$ the Dirac equation formally reduces to a Schr\"odinger equation; the standard derivation splits the four-spinor into ``upper-large'' and ``lower-small'' two-spinors, $\psi_l = e^{\ri c^2 t}(\psi_1,\psi_2)^\top$ and $\psi_s = e^{-\ri c^2 t}(\psi_3,\psi_4)^\top$, and an order-of-magnitude argument in~$c$ eliminates~$\psi_s$ to leave a Schr\"odinger equation for~$\psi_l$~\cite{CLY}. There are thus two limits in play, $\veps\to0$ and $c\to\infty$, and two corresponding ways to arrive at a frozen Gaussian description of the non-relativistic semiclassical particle: one may take the non-relativistic limit of the Dirac equation first and then apply the FGA, or apply the FGA to the Dirac equation first and then take the non-relativistic limit of the resulting representation. The question is whether these two routes agree, whether the diagram of limits commutes.

\begin{thm}[Non-relativistic limit of the FGA \cite{CLY}]
\label{thm:nonrel}
The non-relativistic limit of the Frozen Gaussian Approximation for the field-free linear Dirac equation is the Frozen Gaussian Approximation for the Schr\"odinger equation. That is, the formal diagram relating the semiclassical and non-relativistic limits commutes: the FGA constructed from the non-relativistic Dirac equation coincides with the non-relativistic limit of the FGA constructed from the Dirac equation.
\end{thm}

The proof proceeds by expanding the Dirac Hamiltonian's derivatives in powers of~$c^{-2}$ (e.g., $\partial_P h(Q,P) = P\bigl(1 - P^2/(2c^2) + O(c^{-4})\bigr)$) and tracking these expansions through the flow~\eqref{eq:fga-flow}, the eigenvector formulas, and the amplitude system~\eqref{eq:dirac-amp}, verifying that each ingredient of the Dirac FGA converges to its Schr\"odinger counterpart. The result is the quantum-to-continuum analogue of the asymptotic-preservation Theorem~\ref{thm:ap}: a statement that the approximation respects a second limit built into the problem. The recurrence of the commuting-diagram structure across two entirely different settings (a neural-network solver for Boltzmann in Section~\ref{subsec:ap}, a semiclassical approximation for Dirac here) again suggests that this structure is not tied to a particular method. Rather, it may reflect a more intrinsic feature of the underlying asymptotic bridge.

\begin{figure}[ht]
  \centering
  \includegraphics[width=0.85\textwidth,trim=70 235 70 130,clip]{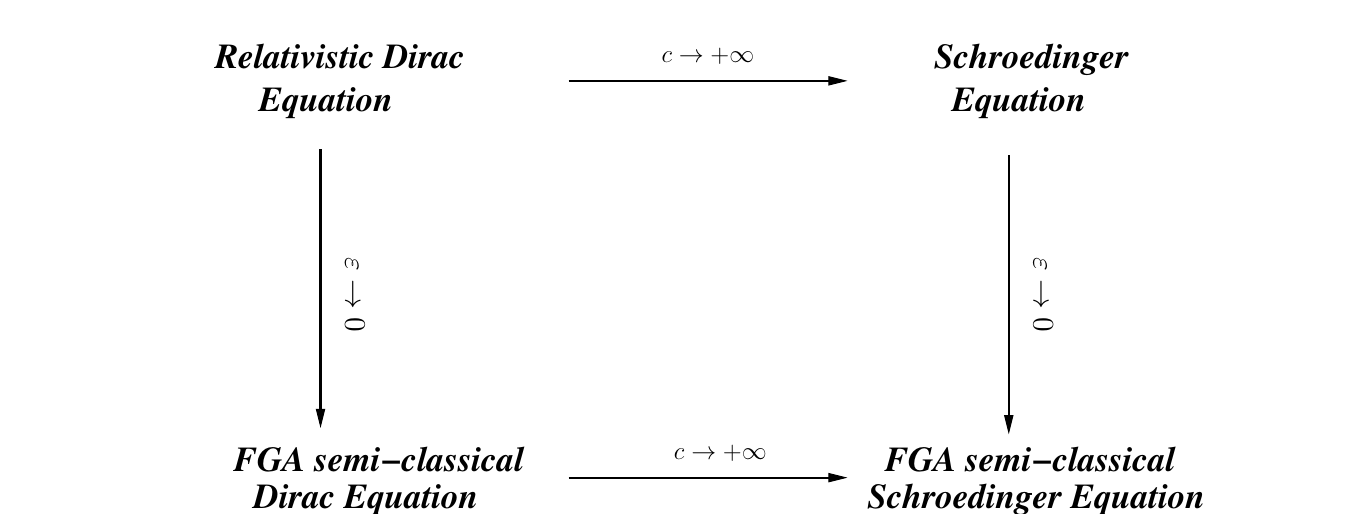}
  \caption{Commutation scheme: the FGA for the Schr\"odinger equation arises both from the semiclassical limit of the non-relativistic Dirac equation and from the non-relativistic limit of the semiclassical (FGA) Dirac description. Reproduced from~\cite{CLY}.}
  \label{fig:commutation}
\end{figure}

\subsection{Error estimate}
\label{subsec:fga-error}
The FGA is not merely a formal construction: it converges, with a rate, to the true solution.

\begin{thm}[FGA convergence \cite{swart2009mathematical, LuYang}]\label{thm:fga-error}
Under suitable regularity and decay hypotheses on the potential and on the initial data, the Frozen Gaussian Approximation satisfies
\begin{equation}\label{eq:fga-rate}
  \bigl\|\psi^\veps(t,\cdot) - \psi^\veps_{\FGA}(t,\cdot)\bigr\|_{L^2}
  \le C\,\veps,
\end{equation}
where the constant~$C$ depends on the dimension, on the regularity of the potential, and on the time horizon~$T$, but not on~$\veps$.
\end{thm}

The estimate~\eqref{eq:fga-rate} holds uniformly through caustics (this is the payoff for which the frozen Gaussian was constructed) and the proof rests on the same oscillatory-integral identities (Lemma~\ref{lem:osc}) that drive the derivation of the amplitude equation, applied now to bound the residual generated by truncating the ansatz at first order.

\subsection{Numerical illustration}
\label{subsec:fga-numerics}
The construction has been implemented and validated for the Dirac equation in~\cite{CLY}. A direct accuracy test compares the FGA solution against a high-accuracy spectral reference solution for a smooth potential; the two are visually indistinguishable at $\veps = 1/2^7$, consistent with the $O(\veps)$ rate of Theorem~\ref{thm:fga-error}. A more demanding test is the \emph{Klein paradox}, the relativistic phenomenon in which a particle incident on a sufficiently steep potential barrier is transmitted rather than reflected, through coupling between the positive- and negative-energy branches. This effect is genuinely quantum: in the semiclassical regime the direct FGA fails to reproduce the transmission, which illustrates a limitation of the method and motivates the further developments discussed in~\cite{CLY}.

\begin{figure}[ht]
  \centering
  \includegraphics[width=0.46\textwidth]{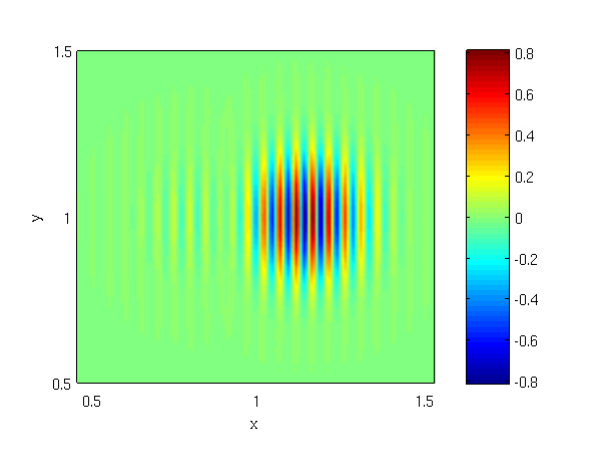}
  \hfill
  \includegraphics[width=0.46\textwidth]{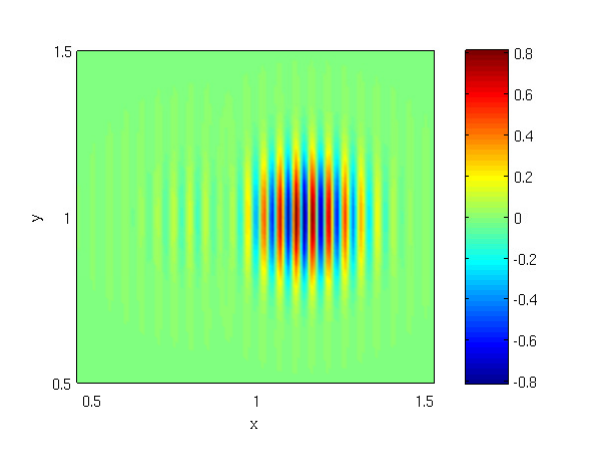}
  \caption{Accuracy test for the FGA applied to the free Dirac equation ($A\equiv 0$, $V\equiv 0$, $c=1$) at $\veps = 1/2^7$, with Gaussian initial spinor of width $\mu = 32$ centered at $x_0 = (1,1,1)$ and momentum $p_0 = (1,0,0)$. The left panel shows the FGA solution and the right panel the reference solution computed by a Fourier spectral method; the two are visually indistinguishable, illustrating the $O(\veps)$ accuracy of Theorem~\ref{thm:fga-error}. Reproduced from~\cite{CLY}.}
  \label{fig:fga-accuracy}
\end{figure}

\begin{figure}[ht]
  \centering
  \includegraphics[width=0.46\textwidth]{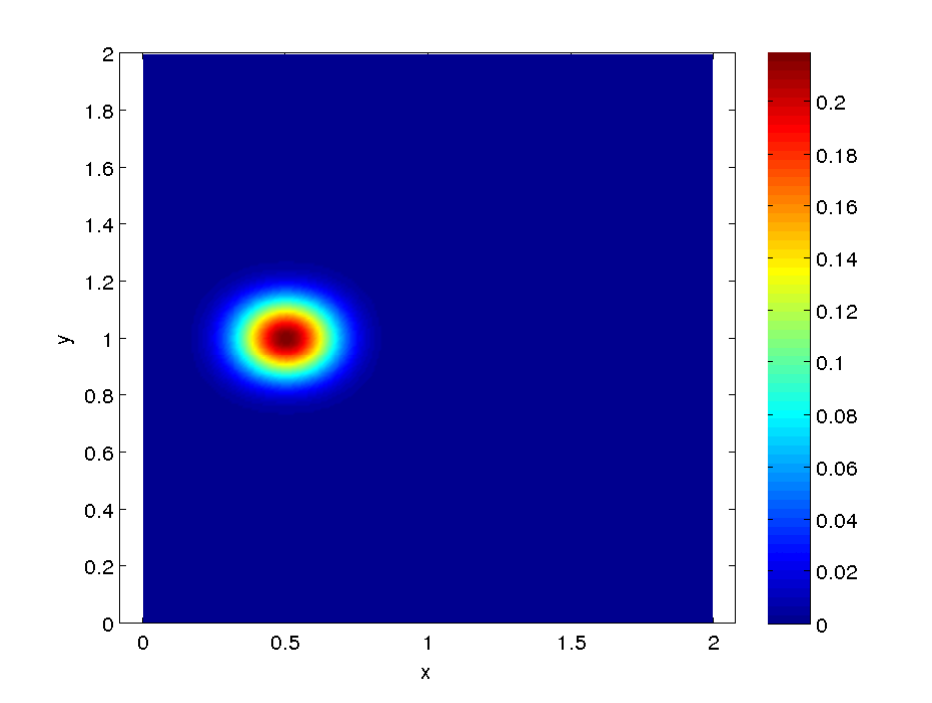}
  \hfill
  \includegraphics[width=0.46\textwidth]{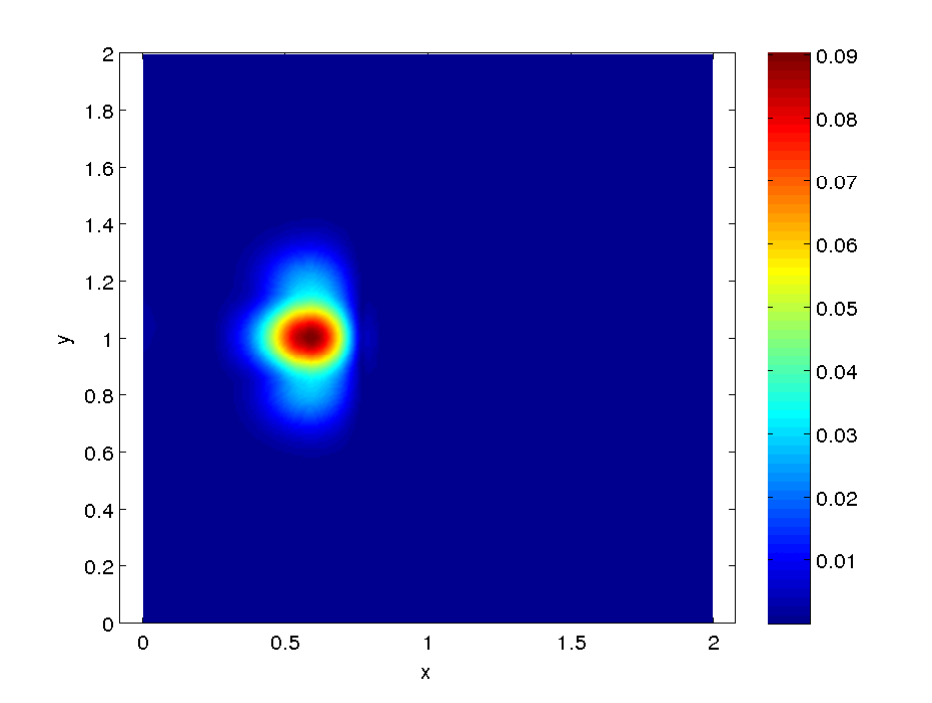}
  \caption{FGA simulation of the Klein paradox in two dimensions in the semiclassical regime ($\veps = 2^{-7}$, $c = 1$), for an incident Gaussian wave packet impinging on a steep barrier potential. The plots show the density $|\bd\psi_1|^2$ of the first spinor component: initial data (left) and the state at $T = 1$ after interaction with the barrier (right). Under the direct FGA the wave packet is almost entirely reflected; the substantial quantum transmission ($\tau \approx 0.485$) predicted in the quantum regime cannot be recovered by a pure semiclassical method, illustrating an intrinsic limitation of the direct FGA in capturing this quantum effect. Reproduced from~\cite{CLY}.}
  \label{fig:klein}
\end{figure}

The non-relativistic limit of Theorem~\ref{thm:nonrel} is also borne out numerically: as the speed of light~$c$ is increased through a range of values, the FGA-for-Dirac solution converges to the FGA-for-Schr\"odinger solution, confirming the commuting diagram of Figure~\ref{fig:commutation}.

\section{Toward AI-assisted asymptotic derivation}\label{sec:ai}
The two bridges of Sections~\ref{sec:kinetic-fluid} and~\ref{sec:fga} were chosen to make a single point vivid: the asymptotic derivation is a structured symbolic object, assembled from a small vocabulary of recurring moves, and the harder its bookkeeping the stronger the case for mechanical assistance. The preceding examples also illustrate what such a mathematical seed should contain: not only final formulas, but also the intermediate identities, projections, side conditions, and limiting diagrams that make the derivation reusable. This final section develops that perspective. We describe what makes such derivations a natural target for AI assistance, sketch a concrete framework \cite{MathEye} for turning imperfect machine-generated proofs into verified ones, and return to the Frozen Gaussian derivation as a candidate for formalization.

\subsection{The bookkeeping problem}
\label{subsec:bookkeeping}
Consider again the derivation of the FGA amplitude equation in Section~\ref{subsec:fga-derivation}. Every step is an instance of a named move: differentiate the phase~\eqref{eq:phase-derivs}; Taylor-expand the potential~\eqref{eq:V-taylor}; apply an oscillatory-integral identity to trade a power of $(x-Q)$ for a power of~$\veps$ (Lemma~\ref{lem:osc}); project onto an eigenspace; collect the surviving order-$\veps$ terms into a transport equation. None of these moves is where the main conceptual novelty lies. Yet the derivation is genuinely error-prone: the holomorphic derivative
$\partial_z = \partial_q - \ri\partial_p$ and the matrix $Z = \partial_z(Q + \ri P)$ must be threaded consistently through every identity; the Einstein summation hides index errors; the Dirac case multiplies the difficulty by carrying a two-dimensional eigenspace whose basis rotates along the flow, so that the matrices $\mathcal{L}$, $\mathcal{M}$, $\mathcal{N}$ in~\eqref{eq:dirac-amp} each require a separate, careful computation. A single dropped term (a side condition for an integration by parts that was never checked, a factor of~$c^{-2}$ misplaced in the non-relativistic expansion) propagates silently through pages of subsequent algebra.

This is the regime in which machine assistance promises the most. The work is symbolic and rule-governed, the rules are few, and the failure mode is loss of bookkeeping rather than loss of insight. The question is how to organize the assistance so that the machine handles the bookkeeping while the human supplies the few genuinely mathematical judgments.

\subsection{Large language models and the verification gap}
\label{subsec:llm-gap}
Large language models can now produce plausible drafts of mathematical arguments, including proof scripts in formal languages such as Lean~\cite{Lean4} and its mathematical library~\cite{Mathlib}, and specialized systems trained for theorem proving (among them DeepSeek-Prover~\cite{DeepSeek} and Llemma~\cite{Llemma}) routinely succeed on competition-style problems. But a draft that looks correct is not a proof. As arguments grow longer and the side conditions subtler, generated scripts fail unpredictably: a missing measurability hypothesis, a misapplied library lemma, an undischarged precondition causes the proof assistant to reject the script, often with an opaque cascade of downstream errors. There is, in short, a gap between the high-level structure that a language model captures well and the low-level rigor that a verifier demands. For an asymptotic derivation, the high-level structure is exactly the sequence of named moves above, and the low-level rigor is exactly the bookkeeping, so the gap is precisely where the difficulty of these derivations lives.

\subsection{Angelic execution and minimal-obligation refinement}
\label{subsec:angelic}
The framework of~\cite{MathEye} is designed to close this gap. Its central idea is \emph{angelic execution}: rather than aborting at the first failed step, the system runs the machine-generated script under an optimistic semantics, replacing each fragment it cannot justify with a symbolic placeholder while preserving the surrounding structure. The result is an \emph{angelic proof}, a structurally coherent script that would be complete if its placeholders were filled, and that exposes exactly which obligations the language model could not discharge.

Presenting all of these obligations to a human would defeat the purpose, since they are typically numerous and interdependent. The framework therefore builds a \emph{proof dependence graph} recording the symbol-level data flow among the fragments, and solves a residual analysis problem (formulated as a MaxSAT optimization) to select a \emph{minimal core} of obligations whose resolution makes all the others follow automatically. The human is queried only about this minimal core, in ordinary mathematical language rather than in the syntax of the proof assistant; once those few questions are answered, the system reconstructs the complete, verified proof. On benchmark problems across analysis and algebra this reduces the number of human interactions sharply while extending the range of theorems that can be carried all the way to a machine-checked conclusion~\cite{MathEye}.

\begin{figure}[ht]
  \centering
  \includegraphics[width=\textwidth]{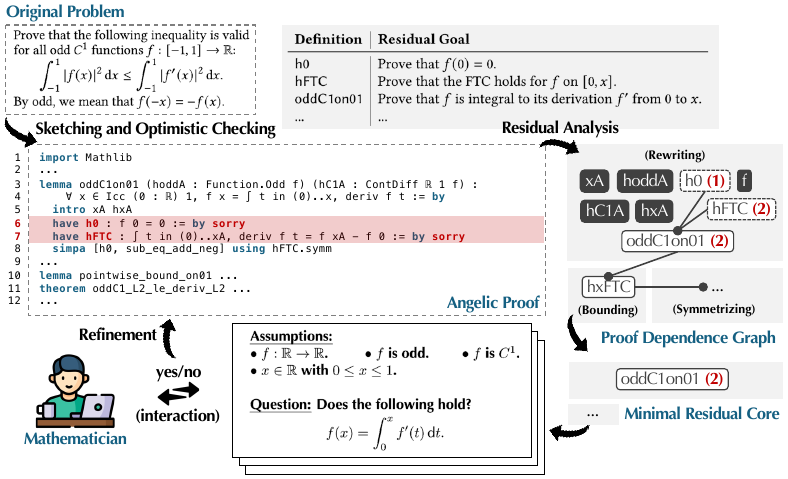}
  \caption{The MathEye workflow for turning a machine-generated proof into a verified one, illustrated on a sample real-analysis problem. The system sketches a proof in Lean and runs it under optimistic (angelic) checking, leaving the steps it cannot yet justify as \texttt{sorry} placeholders (the \emph{angelic proof}). Residual analysis collects these as residual goals and records their dependencies in a proof dependence graph, from which a \emph{minimal residual core} of obligations is selected. The mathematician is then queried only about this core, in ordinary mathematical language rather than proof-assistant syntax; the yes/no answers drive a refinement loop that reconstructs the complete, verified proof. Reproduced from~\cite{MathEye}.}
  \label{fig:angelic}
\end{figure}

The relevance to asymptotic derivation is direct. The named moves of Section~\ref{subsec:bookkeeping} are the high-level structure a language model can be expected to propose; the unchecked side conditions and index manipulations are the placeholders an angelic execution would expose; and the few genuinely mathematical judgments (does this integration by parts have the regularity it needs, is this the correct eigenprojection) are the minimal core a mathematician should be asked to confirm.

\subsection{The reusable-formalization principle}
\label{subsec:transfer}
We can now state this perspective in its sharpest form. Asymptotic derivations across the kinetic-to-fluid and quantum-to-continuum bridges (e.g., in Sections~\ref{subsec:ap} and~\ref{subsec:nonrel}) are assembled from the same small vocabulary: power matching, oscillatory-integral identities, integration by parts, the Fredholm alternative, moment closure, eigenprojection, and changes of variables along characteristic flows. Once a single nontrivial derivation in this vocabulary has been formalized (its lemmas stated and certified, its dependency structure recorded, its recurring side conditions catalogued) the formalization of the next derivation in the same family reuses most of that infrastructure. The certified oscillatory-integral identity that powers the FGA amplitude equation is the same identity, up to the nouns, that one needs for the Gaussian beam method, for the WKB expansion, and for the general linear hyperbolic FGA. Formalizing a difficult member of a family (e.g., the Dirac FGA) can make its easier relatives substantially less costly to certify. We call this the \emph{reusable-formalization principle}: formalizing a difficult member of a family can lower the marginal cost of certifying related derivations, making the first hard derivation the one most worth investing in.

\subsection{A test case: one oscillatory-integral identity}
\label{subsec:primitive}
To make the proposal concrete, consider what it would take to certify a single primitive: the first oscillatory-integral identity of Lemma~\ref{lem:osc}, the move that trades a factor of $(x-Q)$ for a factor of~$\veps$. Writing only its \emph{statement} in the language of a proof assistant already does useful work, because it forces every side condition that a hand derivation leaves implicit to be named. Schematically, in a Lean-like syntax (ASCII for the symbols):

{\footnotesize
\begin{verbatim}
-- Lemma 3.2, identity (1): trade (x-Q) for a power of e.
-- OscEquiv e Phi f g  :  the integrals of f and of g against
--   exp(i*Phi/e) agree.    dz := d/dq - i * d/dp.
-- x and t are fixed parameters; dotBilinear has no conjugation.
theorem osc_linear
    (he   : 0 < e)
    (hPhi : IsFGAPhase Phi Q P S) -- flow and action identities
    (hFlow: UniformFlowBounds Q P)
    (hZ   : InvertibleWithBddDerivs Z)
    (hw   : BddDerivs w)          -- amplitude class, Lemma 3.2
    : OscEquiv e Phi
        (fun z => dotBilinear (w z) (x - Q z))
        (fun z => -e * divz
          (fun s => transpose (inv (Z s)) * w s) z)
\end{verbatim}
}

These schematic hypotheses record the Hamiltonian phase and action identities, the uniform bounds on the flow derivatives and on~$Z^{-1}$, and the amplitude regularity from Lemma~\ref{lem:osc}. Here \texttt{divz} denotes divergence with respect to~$\partial_z$, and the transpose implements the index contraction $\sum_j w_jZ^{-1}_{jk}$ in~\eqref{eq:osc-linear}. Integration by parts is justified first on Schwartz data using suitable cutoffs, followed by the oscillatory-integral limit and the $L^2$ bounds. Stating the lemma in this form makes these analytic obligations explicit; the displayed syntax is schematic rather than an executable Lean theorem.

The proof is another matter. A proof assistant such as Lean~\cite{Lean4} supplies, through Mathlib~\cite{Mathlib}, general analytic tools (e.g., Schwartz space and integration by parts), but not yet a developed theory of oscillatory integrals with complex Gaussian phase, nor the
stationary-phase and non-degeneracy estimates on which the proof of Lemma~\ref{lem:osc} actually rests; certifying this one identity would mean building that layer. This is precisely the kind of initial formalization effort described by the reusable-formalization principle. Once the oscillatory-integral machinery has been built, it becomes a reusable component for many derivations in the same family. Thus the first difficult primitive is also the one whose certification is most valuable.

\subsection{Outlook: an AI-assisted FGA derivation}
\label{subsec:fga-ai}
Concretely, one can imagine an AI-assisted derivation of the amplitude equation~\eqref{eq:fga-amp-schr}, and of its Dirac counterpart~\eqref{eq:dirac-amp}, proceeding along the following lines. The named moves of Section~\ref{subsec:bookkeeping} could be encoded as a library of certified primitives in a proof assistant: the phase-derivative identities~\eqref{eq:phase-derivs}, the Taylor expansion with explicit remainder, and above all the oscillatory-integral identities of Lemma~\ref{lem:osc}, with their hypotheses made explicit. A language model proposes a derivation sketch (the sequence of moves, in order), which angelic execution renders as an angelic proof, isolating as placeholders the steps whose side conditions it cannot discharge: typically the regularity needed for each integration by parts and the eigenprojection within the degenerate Dirac eigenspace. Residual analysis then extracts the minimal core of these placeholders, the mathematician confirms them, and the system reconstructs a verified derivation of the transport equation. The certified library that results would be reusable: it would be much the same library that could carry a formalized Gaussian beam derivation, or a formalized WKB expansion, or the convergence estimate of Theorem~\ref{thm:fga-error}, whose proof rests on the very same identities. The precise architecture of such a library remains open.
 
From this viewpoint, the FGA derivation is not only an example to be checked, but also a starting point for a future library of semiclassical asymptotic arguments. Its phase identities, Taylor expansions, oscillatory-integral lemmas, and projection steps are precisely the kinds of components that can be reused in related derivations.
 
This is one possible direction suggested by the examples above. The claim is not that asymptotic analysis will be automated away, but the opposite: that the discipline becomes \emph{more} robust, not less, as these tools mature, because each derivation that is formalized becomes verified, reusable infrastructure rather than a fragile pencil-and-paper artifact to be re-checked by hand at every reuse.

\section{Conclusion and outlook}\label{sec:conclusion}
A single methodological pattern has run through this paper. The passage from kinetic to fluid, reviewed through the diffusion limit of radiative transfer and its $H$-function interface conditions; the analysis of physics-informed neural networks for Boltzmann equations, with its micro--macro decomposition and asymptotic-preserving estimates; the passage from quantum to continuum, developed in depth through the Frozen Gaussian Approximation, its Dirac extension, and its non-relativistic limit: all are instances of the same recipe of ansatz, substitution, matching, and assembly, drawn from the same small vocabulary of symbolic moves.

That vocabulary is the paper's bridge to a second theme. Because the moves are few, rule-governed, and error-prone in ways that are well suited to machine assistance, the asymptotic derivation is a natural target for formalization. Moreover, because the moves recur across derivations, a single formalized derivation yields reusable infrastructure that lowers the cost of all the others: the reusable-formalization principle of Section~\ref{subsec:transfer}.

The longer-term vision is not the replacement of the analyst, but the strengthening of analysis. As large language models, proof assistants, and frameworks such as angelic execution~\cite{MathEye} mature, carefully written derivations may become permanent, composable, and reusable mathematical objects. In this way, review articles may also take on a new role: they can plant structured mathematical seeds in the emerging AI-assisted ecosystem, making the bridges between scales not only understood, but increasingly certified and easier to extend.

\section*{Acknowledgements}
The author warmly thanks his collaborators on the works reviewed here (in particular Elie Abdo, Lihui Chai, Carlos Garc\'ia-Cervera, Fran\c{c}ois Golse, Ruimeng Hu, Zhongyi Huang, Shi Jin, Xiaomei Liao, Emmanuel Lorin, Jianfeng Lu, and Hao Wu) and the organizers of ICCM 2025 for the invitation to deliver the lecture on which this paper is based. He also thanks his collaborators on the MathEye framework: Yanju Chen, Ruizhe Qian, Jiaming Shan, Yufei Ding, Osbert Bastani, and Yu Feng.

\bibliographystyle{amsplain}
\bibliography{yang-iccmp2025}

\end{document}